# An Infinite-Time Relaxation Theorem for Differential Inclusions

Brian Ingalls and Eduardo D. Sontag[*]
Dept. of Mathematics, Rutgers University, NJ
{ingalls,sontag}@math.rutgers.edu

Yuan Wang[†]
Dept. of Mathematics, Florida Atlantic University, FL
ywang@math.fau.edu

November 13, 2018

## Abstract

The fundamental relaxation result for Lipschitz differential inclusions is the Filippov-Ważewski Relaxation Theorem, which provides approximations of trajectories of a relaxed inclusion on finite intervals. A complementary result is presented, which provides approximations on infinite intervals, but does not guarantee that the approximation and the reference trajectory satisfy the same initial condition.

## 1 Introduction

This note studies the approximation of solutions of the relaxation of a differential inclusion of the type:

$$\dot{x}(t) \in F(t, x(t)) \tag{1}$$

where the set-valued function $F$ is locally Lipschitz and takes values which are nonempty and closed. The relaxation considered is the inclusion

$$\dot{x}(t) \in \text{clco } F(t, x(t)) \tag{2}$$

where clco stands for closed convex hull.

The fundamental result on approximations of solutions of (2) by solutions of (1) is the Filippov-Ważewski Relaxation Theorem (cf. [2, 3, 5, 6, 7]). This result says that the solution set of (1) is dense in the solution set of (2) in the topology of uniform convergence on compact intervals. (The paper [7] provides a continuous version of the Theorem which is closely related to the tools used in this note).

In particular, the Filippov-Ważewski Theorem says that given a trajectory of the relaxed system (2) defined on a finite interval, there exists a trajectory of (1) *with the same initial condition* which approximates the trajectory of the relaxed system *on that finite interval*. A complementary result is presented in this note. Roughly speaking, it is shown that the solution set of initial value problems of the type

$$\dot{x}(t) \in F(t, x(t)), \qquad x(0) = \xi_1, \tag{3}$$

is dense in the solution set of initial value problems of the type

$$\dot{x}(t) \in \text{ clco } F(t, x(t)), \qquad x(0) = \xi_2, \tag{4}$$

in the "$C^0$ Whitney topology" on the *infinite* interval $[0, \infty)$. This is not a generalization of the Filippov-Ważewski theorem, as, given a trajectory of (4), this result does not guarantee the existence of an approximating trajectory of (3) with $\xi_1 = \xi_2$, but rather only with $\xi_1$ arbitrarily close to $\xi_2$.

[*]Supported in part by US Air Force Grant F49620-98-1-0242
[†]This research was supported in part by NSF Grant DMS-9457826

The result in this note provides the existence of trajectories which are approximations in weighted norms on $[0, \infty)$, for example $|f| := \sup_{t \geq 0}\{|f(t)| e^t\}$. Indeed, given any $r : \mathbb{R}_{\geq 0} \to \mathbb{R}_{>0}$, there is an approximation in the norm $|f| := \sup_{t \geq 0}\{|f(t)| r(t)\}$. This is achieved by demanding that the approximation lie in a tube around the reference trajectory which has possibly vanishing radius. An immediate corollary is that the relaxation (2) is forward complete if and only if the inclusion (1) is forward complete.

Included also is a counter-example which shows that one cannot achieve an approximation on the infinite interval if one insists that the approximation satisfy the same initial condition as the reference trajectory.

The motivation for this work was a question in the stability of differential inclusions. It was shown in [8] that a differential inclusion $\dot{x} \in F(x)$ is globally asymptotically stable if and only if it is uniformly globally asymptotically stable, provided that the set-valued map $F$ admits a parameterization of the form $F(x) = \{f(x, u) \; : \; u \in U\}$ where $f(\cdot, \cdot)$ is locally Lipschitz and $U$ is compact. (See [1] for a more general result). The proof in this note combines the tools used in [8] with the main result in the excellent paper [4] which provides continuous selections of solutions of (1).

## 1.1 Basic Definitions and Notation

For each $T > 0$, let $\mathcal{L}[0, T]$ be the $\sigma$-field of Lebesgue measurable subsets of $[0, T]$. Let $X$ be a separable Banach space, whose norm is denoted simply by $|\cdot|$. Let $\mathcal{P}(X)$ denote the family of all nonempty closed subsets of $X$. We use $\mathcal{B}(X)$ for the family of Borel subsets of $X$.

For each interval $\mathcal{I} \subseteq [0, \infty)$, let $L^1(\mathcal{I}, X)$ be the Banach space of Bochner integrable functions $u : \mathcal{I} \to X$ with norm $\|u\| = \int_{\mathcal{I}} |u(t)| \, dt$, and let $L^1_{\text{loc}}(\mathcal{I}, X)$ be the corresponding space of locally integrable functions. Let $AC(\mathcal{I}, X)$ be the Banach space of absolutely continuous functions $u : \mathcal{I} \to X$ with the norm $\|u\|_{AC} = |u(0)| + \|\dot{u}\|$.

We define the distance from a point $\xi \in X$ to a set $K \in \mathcal{P}(X)$ as

$$d(\xi, K) := \inf\{|\xi - \eta| \; : \; \eta \in K\}.$$

For a set $A \in \mathcal{P}(X)$, let $B(A, r)$ denote the set $\{\xi \in X \; : \; d(\xi, A) \leq r\}$. For singleton $A = \{\xi\}$ we write $B(\xi, r)$. For each set $A$ and each constant $c \in \mathbb{R}$, we denote $cA = \{c\xi \; : \; \xi \in A\}$.

**Definition 1.1** The *Hausdorff distance* between two sets $K, L \in \mathcal{P}(X)$ is defined as

$$d_H(K, L) := \max\left\{\sup_{\xi \in K} d(\xi, L), \sup_{\eta \in L} d(\eta, K)\right\}.$$

□

**Definition 1.2** Let $\mathcal{O}$ be an open subset of $X$. Let $\mathcal{I} \subseteq \mathbb{R}_{\geq 0}$ be an interval. The set-valued map $F : \mathcal{I} \times X \to \mathcal{P}(X)$ is said to be *locally Lipschitz* on $\mathcal{O}$ if, for each $\xi \in \mathcal{O}$, there exists a neighbourhood $U \subset \mathcal{O}$ of $\xi$ and a $k_U \in L^1(\mathcal{I}, \mathbb{R})$ so that for any $\eta$, $\zeta$ in $U$,

$$d_H(F(t, \eta), F(t, \zeta)) \leq k_U(t) |\eta - \zeta| \qquad \text{a.e. } t \in \mathcal{I}.$$

$\square$

**Definition 1.3** Let $\mathcal{I} \subseteq [0, \infty)$ be an interval. A function $x : \mathcal{I} \to X$ is said to be a *solution of the differential inclusion*

$$\dot{x}(t) \in F(t, x(t)) \tag{5}$$

if it is absolutely continuous and satisfies (5) for almost every $t \in \mathcal{I}$.

For $T > 0$, a function $x : [0, T) \to X$ is called a *maximal solution of the differential inclusion* if it does not have an extension which is a solution in $X$. That is, either $T = \infty$ or there does not exist a solution $y : [0, T_+] \to X$ with $T_+ > T$ so that $y(t) = x(t)$ for all $t \in [0, T)$. $\square$

**Definition 1.4** A differential inclusion is called *forward complete* if every maximal solution is defined on the interval $[0, \infty)$. $\square$

## 2 Continuous Selections of Trajectories

We begin by presenting a particular case of the main theorem in [4]. In the spirit of keeping this work self-contained, the full statement of the theorem in [4] is included in the appendix.

**Lemma 2.1** Let $T > 0$ and a set-valued map $F : [0, T] \times X \to \mathcal{P}(X)$ be given, and consider the initial value problems for $\eta \in X$,

$$\dot{x} \in F(t, x), \quad x(0) = \eta, \qquad \text{for } t \in [0, T]. \tag{6}$$

Fix $\xi_0 \in X$, and suppose $\overline{y} : [0, T] \to X$ is a solution of (6) with $\eta = \xi_0$. Then, if $F$ satisfies

(H1) $F$ is $\mathcal{L}[0, T] \otimes \mathcal{B}(X)$ measurable;

(H2) there exists $k \in L^1([0, T], \mathbb{R})$ such that for any $\xi, \eta \in X$,

$$d_H(F(t, \xi), F(t, \eta)) \leq k(t) |\xi - \eta| \qquad \text{a.e. } t \in [0, T];$$

(H3) there exists a point $x_0 \in X$ and a $\beta_0 \in L^1([0, T], \mathbb{R})$ such that

$$d(x_0, F(t, x_0)) \leq \beta_0(t) \qquad \text{a.e. } t \in [0, T];$$

it follows that for each $\varepsilon_0 > 0$ there exists a function $x : [0, T] \times X \to X$ such that

(a) For every $\eta \in X$, the function $t \mapsto x(t, \eta)$ is a solution of (6);

(b) the map $\eta \mapsto x(\cdot, \eta)$ is continuous from $X$ into $AC([0, T], X)$;

(c) for each $\eta \in X$, and each $t \in [0, T]$,

$$|\overline{y}(t) - x(t, \eta)| \leq (\varepsilon_0 + |\xi_0 - \eta|) e^{\int_0^t k(s)\,ds}.$$

*Proof.* We apply Theorem 3.1 of [4] (which is included in the appendix) with $S = X$. Define a continuous map $\eta \mapsto y(\cdot, \eta)$ from $X$ into $AC([0, T], X)$ by the constant assignment $\eta \mapsto \overline{y}(\cdot)$. The Theorem requires a continuous map $\beta_y : X \to L^1([0, T], \mathbb{R})$ so that for each $\eta \in X$,

$$d(\dot{y}(t, \eta), F(t, y(t, \eta))) \leq \beta_y(\eta)(t) \qquad \text{a.e.} \quad t \in [0, T].$$

Since we have $y(t, \eta) \equiv \overline{y}(t)$ for each $\eta$, and $\overline{y}(\cdot)$ is a solution of (6), it follows that we may choose $\beta_y(\eta) = 0$. We note also that our hypothesis (H3) is equivalent to the hypothesis (H4$_0$) in [4], as the global Lipschitz condition (H2) gives

$$d(0, F(t, 0)) \leq |x_0| + k(t) |x_0| + d(x_0, F(t, x_0)) \qquad \text{a.e.} \quad t \in [0, T]$$

for any $x_0 \in X$. Then, for any given $\varepsilon_0 > 0$, the Theorem provides the existence of a function $x : [0, T] \times X \to X$ which satisfies (a) and (b) above, as well as

$$|(\overline{y}(t) - x(t, \eta)) - (\xi_0 - \eta)| \leq (\varepsilon_0 + |\xi_0 - \eta|) e^{\int_0^t k(s)\,ds} - |\xi_0 - \eta|,$$

for each $\eta \in X$ and each $t \in [0, T]$, from which (c) follows easily. ∎

## 3 Approximations of Trajectories of Relaxed Inclusions

We next state a lemma on continuous selections of approximations of a trajectory of a relaxed inclusion on a finite interval.

**Lemma 3.1** Let $T > 0$, $\xi_0 \in X$, and a set-valued map $F : [0, T] \times X \to \mathcal{P}(X)$ be given, and consider the initial value problems for $t \in [0, T]$

$$\dot{x} \in F(t, x), \quad x(0) = \xi_0, \tag{7}$$

and

$$\dot{x} \in \text{clco } F(t, x), \quad x(0) = \xi_0. \tag{8}$$

Suppose $z : [0, T] \to X$ is a solution of (8), and let $\varepsilon > 0$ be given. Let

$$\mathcal{T} := \{\xi \in X \ : \ |\xi - z(t)| \leq \varepsilon \ \text{ for some } \ t \in [0, T]\},$$

the $\varepsilon$-tube around the image of $z$. Then, if $F$ satisfies

(H1) $F$ is $\mathcal{L}[0, T] \otimes \mathcal{B}(X)$ measurable;

(H2′) there exists $k_0 \in L^1([0, T], \mathbb{R})$ such that for any $\xi, \eta \in B(\mathcal{T}, 1)$

$$d_H(F(t, \xi), F(t, \eta)) \leq k_0(t) |\xi - \eta| \qquad \text{a.e.} \ t \in [0, T];$$

(H3′) there exists $\alpha \in L^1([0, T], \mathbb{R})$ such that for each $\xi \in B(\mathcal{T}, 1)$

$$\sup\{|\zeta| \ : \ \zeta \in F(t, \xi)\} \leq \alpha(t) \qquad \text{a.e.} \ t \in [0, T];$$

it follows that there exists a $\delta > 0$ and a function $x : [0,T] \times V \to X$, where $V := B(\xi_0, \delta)$ such that

(a) For every $\eta \in V$, the function $t \mapsto x(t, \eta)$ is a solution of the initial value problem

$$\dot{x} \in F(t, x), \quad x(0) = \eta, \quad \text{for } t \in [0, T]; \tag{9}$$

(b) the map $\eta \mapsto x(\cdot, \eta)$ is continuous from $V$ into $AC([0,T], X)$;

(c) for each $\eta \in V$,

$$|z(t) - x(t, \eta)| \leq \varepsilon \quad \forall t \in [0, T].$$

*Proof.* We combine Lemma 2.1 with the Filippov-Ważewski Relaxation Theorem (for the statement of the Relaxation Theorem in the full generality used here, see e.g. [7]). Let $T$, $F$, $\xi_0$ and $z$ be as above, and let $\varepsilon > 0$ be given. By the Relaxation Theorem, there exists a solution $\bar{y}$ of (7) which satisfies

$$|z(t) - \bar{y}(t)| < \frac{\varepsilon}{2} \quad \forall t \in [0, T].$$

Next we turn to Lemma 2.1. To apply the Lemma, we need to modify the function $F$ to ensure the Lipschitz and boundedness properties hold over the whole space $X$.

We define $\Phi : X \to [0, 1]$ by

$$\Phi(x) := \max\{1 - d(x, \mathcal{T}), 0\}.$$

Set $\widetilde{F}(t, x) := \Phi(x) F(t, x)$. Since $\widetilde{F}(t, x)$ and $F(t, x)$ agree for any $(t, x) \in [0, T] \times \mathcal{T}$, it follows that the trajectories of (7), (8), and (9) which stay inside $\mathcal{T}$ are the same as those of the differential inclusions with $\widetilde{F}(t, x)$ in the place of $F(t, x)$. Moreover, $\widetilde{F}$ satisfies the hypotheses (H1)–(H3) of Lemma 2.1 as follows: (H1) is immediate. (H3) follows by taking any $x_0 \in \mathcal{T}$ and choosing $\beta_0(t) = |x_0| + \alpha(t)$. For (H2), we find, for any $t \in [0, T]$, for each pair $x, y \in X$,

i) if $x, y \in B(\mathcal{T}, 1)$,

$$\begin{aligned}
d_H(\widetilde{F}(t,x), \widetilde{F}(t,y)) &= d_H(\Phi(x)F(t,x), \Phi(y)F(t,y)) \\
&\leq d_H(\Phi(x)F(t,x), \Phi(y)F(t,x)) + d_H(\Phi(y)F(t,x), \Phi(y)F(t,y)) \\
&\leq |\Phi(x) - \Phi(y)| \sup\{|\zeta| : \zeta \in F(t,x)\} + |\Phi(y)| d_H(F(t,x), F(t,y)) \\
&\leq |x - y|\alpha(t) + d_H(F(t,x), F(t,y)) \\
&\leq (\alpha(t) + k_0(t)) |x - y|;
\end{aligned}$$

ii) if $x \in B(\mathcal{T}, 1)$, $y \notin B(\mathcal{T}, 1)$,

$$\begin{aligned}
d_H(\widetilde{F}(t,x), \widetilde{F}(t,y)) &= d_H(\Phi(x)F(t,x), \{0\}) \\
&= \sup\{|\zeta| : \zeta \in \Phi(x)F(t,x)\} \\
&= |\Phi(x)| \sup\{|\zeta| : \zeta \in F(t,x)\} \\
&\leq |\Phi(x)| \alpha(t) \\
&= |\Phi(x) - \Phi(y)| \alpha(t) \\
&\leq |x - y| \alpha(t);
\end{aligned}$$

iii) if $x, y \notin B(\mathcal{T}, 1)$, $d_H(\widetilde{F}(t,x), \widetilde{F}(t,y)) = d_H(\{0\}, \{0\}) = 0$.

Hence the global Lipschitz condition (H2) holds with $k(t) = \alpha(t) + k_0(t)$.

We apply the Lemma with $\overline{y}$ and $\varepsilon_0 := \frac{\varepsilon}{4m}$, where $m := \exp(\int_0^T k(s)\,ds)$. Since its image lies in $\mathcal{T}$, $\overline{y}$ is a trajectory of (7) with $\widetilde{F}(t,x)$ in the place of $F(t,x)$. The Lemma gives the existence of a function $x : [0,T] \times X \to X$ so that

(a) For every $\eta \in X$, the function $t \mapsto x(t, \eta)$ is a solution of
$$\dot{x} \in \widetilde{F}(t, x), \quad x(0) = \eta; \tag{10}$$

(b) the map $\eta \mapsto x(\cdot, \eta)$ is continuous from $X$ into $AC([0,T], X)$;

(c) for each $\eta \in X$,
$$|\overline{y}(t) - x(t, \eta)| \leq (\varepsilon_0 + |\xi_0 - \eta|) e^{\int_0^t k(s)\,ds} \quad \forall t \in [0, T].$$

Choosing $\delta = \frac{\varepsilon}{4m}$, we have from (c) that for each $\eta \in V := B(\xi_0, \delta)$,
$$|\overline{y}(t) - x(t, \eta)| \leq \frac{\varepsilon}{2} \quad \forall t \in [0, T].$$

Thus for each $\eta \in V$,
$$|z(t) - x(t, \eta)| \leq |z(t) - \overline{y}(t)| + |\overline{y}(t) - x(t, \eta)| < \varepsilon \quad \forall t \in [0, T].$$

This implies that for each $\eta \in V$, the trajectory $x(\cdot, \eta)$, lies in the tube $\mathcal{T}$ in which $\widetilde{F}$ and $F$ coincide, so these are in fact trajectories of the original system. We conclude that the restriction of $x$ to $[0, T] \times V$ satisfies the required conditions. ∎

Our main result will be an immediate corollary of the following technical lemma. Given $0 < T \leq \infty$ and a trajectory $z : [0, T) \to X$ of the relaxed system (2), this lemma will show the existence of, for any strictly increasing sequence of times $T_k \to T$ and for each nonnegative integer $k$, a sequence $\{\eta_j^k\}_{j=1}^\infty$, whose elements are close to $z(T_k)$ and which satisfy a continuous reachability property.

**Lemma 3.2** Let $0 < T \leq \infty$, and suppose the set-valued map $F : [0, T) \times X \to \mathcal{P}(X)$ satisfies the following properties.

(H1″) $F$ is $\mathcal{L}[0, T) \otimes \mathcal{B}(X)$ measurable;

(H2″) for each $R > 0$, there exists $k_R \in L^1_{\text{loc}}([0, T), \mathbb{R})$ such that for any $\xi, \eta \in B(0, R)$
$$d_H(F(t, \xi), F(t, \eta)) \leq k_R(t) |\xi - \eta| \quad \text{a.e. } t \in [0, T);$$

(H3″) for each $R > 0$, there exists $\alpha_R \in L^1_{\text{loc}}([0, T), \mathbb{R})$ such that for each $\xi \in B(0, R)$
$$\sup\{|\zeta| : \zeta \in F(t, \xi)\} \leq \alpha_R(t) \quad \text{a.e. } t \in [0, T);$$

Fix $\xi \in X$ and let $z : [0, T) \to X$ be a solution of
$$\dot{x} \in \text{clco } F(t, x), \quad x(0) = \xi.$$

Let $r : [0, T) \to \mathbb{R}$ be a continuous function satisfying $r(t) > 0$ for all $t \in [0, T)$. Let $\{T_k\}_{k=0}^\infty$ be any strictly increasing sequence of times so that $T_0 = 0$ and $T_k \to T$ as $k \to \infty$. Then there exists a sequence $\{\delta_k\}_{k=0}^\infty$ of positive numbers and, for each nonnegative integer $k$, a sequence of points $\{\eta_j^k\}_{j=1}^\infty$ which satisfy the following.

- for each $k \geq 0$, $\delta_k \leq \min\{r(t) : t \in [T_k, T_{k+1}]\}$;

- for each $k \geq 0$, $\eta_j^k \in V_k := B(z(T_k), \delta_k)$ for all $j \geq 1$;

- for any $k \geq 1$, *if* a subsequence $\{\eta_{j_l}^k\}_{l=1}^\infty$ converges, say to $\eta^k$, then the subsequence $\{\eta_{j_l}^{k-1}\}_{l=1}^\infty$ also converges, say to $\eta^{k-1}$, and there is a solution $\bar{x} : [0, T_k - T_{k-1}] \to X$ of the initial value problem

$$\dot{x} \in F(T_{k-1} + t, x), \quad x(0) = \eta^{k-1} \quad \text{for } t \in [0, T_k - T_{k-1}], \tag{11}$$

which satisfies

$$|\bar{x}(t) - z(T_{k-1} + t)| \leq r(T_{k-1} + t) \quad \forall t \in [0, T_k - T_{k-1}], \tag{12}$$

and has $\bar{x}(T_k - T_{k-1}) = \eta^k$.

*Proof.* For each positive integer $k$, let

$$r_k := \min\{r(t) : t \in [T_{k-1}, T_k]\}.$$

For each positive integer $k$ we will build a family of trajectories which approximate $z$ on the time interval $[T_{k-1}, T_k]$. On each such interval, we will consider the differential inclusions in *backward time*.

We will apply Lemma 3.1 to the problems

$$\dot{x} \in -F(T_k - t, x), \quad \text{for } t \in [0, T_k - T_{k-1}], \tag{13}$$

and

$$\dot{x} \in \text{clco } -F(T_k - t, x), \quad \text{for } t \in [0, T_k - T_{k-1}], \tag{14}$$

for each $k \geq 1$ with appropriate initial conditions.

Set $\delta_0 = r_1$ and $V_0 := B(z(0), r_1)$. We will construct, by induction, for each positive integer $k$,

- a $\delta_k > 0$ which satisfies $\delta_k \leq r_{k+1}$,

- a set $V_k := B(z(T_k), \delta_k)$, and

- a function $x_k : [0, T_k - T_{k-1}] \times V_k$ which satisfies the following:

    (a) for every $\eta \in V_k$, the function $t \mapsto x_k(t, \eta)$ is a solution of (13) with initial condition $x(0) = \eta$;

    (b) the map $\eta \mapsto x_k(\cdot, \eta)$ is continuous from $V_k$ into $AC([0, T_k - T_{k-1}], X)$;

    (c) for each $\eta \in V_k$,

    $$|z(T_k - t) - x_k(t, \eta)| \leq r_k \quad \forall t \in [0, T_k - T_{k-1}].$$

    (d) for every $\eta \in V_k$,

    $$x_k(T_k - T_{k-1}, \eta) \in V_{k-1}.$$

We first make the construction for $k = 1$. Note that, by definition, $z(T_1 - t)$ is a solution of (14) with initial value $z(T_1)$. Further, the hypotheses of Lemma 3.1 are satisfied by the function $F$, and hence $-F$, since there exists $R$ large enough so that $B(0, R)$ contains the image of $z(T_1 - t)$ over $t \in [0, T_1]$. Applying the Lemma with $\varepsilon = r_1$, it follows that there exists a $\delta_1 > 0$ and a function $x_1 : [0, T_1] \times V_1$, where $V_1 := B(z(T_1), \delta_1)$, which satisfies

(a) For every $\eta \in V_1$, the function $t \mapsto x_1(t, \eta)$ is a solution of
$$\dot{x} \in -F(T_1 - t, x), \qquad x(0) = \eta;$$

(b) the map $\eta \mapsto x_1(\cdot, \eta)$ is continuous from $V_1$ into $AC([0, T_1], X)$;

(c) for each $\eta \in V_1$,
$$|z(T_1 - t) - x_1(t, \eta)| \leq r_1 \qquad \forall t \in [0, T_1],$$

(d) for each $\eta \in V_1$,
$$x_1(T_1, \eta) \in V_0.$$

where (c) follows from the choice of $\varepsilon = r_1$, and (d) follows from evaluating (c) at $t = T_1$.

Now, supposing that for some $k \geq 1$ there exist $\delta_k$ and $x_k$ as above, we produce $\delta_{k+1}$ and $x_{k+1}$ as follows.

Consider the function $z(T_{k+1} - t)$ on the interval $t \in [0, T_{k+1} - T_k]$. This solves (14) (for $k+1$) with initial value $z(T_{k+1})$. We apply Lemma 3.1, with $\varepsilon = \min\{\delta_k, r_{k+1}\}$, to find a $\delta_{k+1} > 0$ and a function $x_{k+1} : [0, T_{k+1} - T_k] \times V_{k+1}$, where $V_{k+1} := B(z(T_{k+1}), \delta_{k+1})$, which satisfies

(a) for every $\eta \in V_{k+1}$, the function $t \mapsto x_{k+1}(t, \eta)$ is a solution of (13) (for $k+1$) with initial condition $x(0) = \eta$;

(b) the map $\eta \mapsto x_{k+1}(\cdot, \eta)$ is continuous from $V_{k+1}$ into $AC([0, T_{k+1} - T_k], X)$;

(c) for each $\eta \in V_{k+1}$,
$$|z(T_{k+1} - t) - x_{k+1}(t, \eta)| \leq \varepsilon \leq r_{k+1} \qquad \forall t \in [0, T_{k+1} - T_k],$$

(d) for each $\eta \in V_{k+1}$,
$$x_{k+1}(T_{k+1} - T_k, \eta) \in V_k.$$

where (d) follows from $\varepsilon \leq \delta_k$. Then, by induction, we conclude that there exist such $\delta_k$ and $x_k$ for each $k \geq 1$.

Next, for each positive integer $k$ we consider the concatenated trajectory $y_k : [0, T_k] \to X$ defined by

$$y_k(t) := \begin{cases} x_k(t, z(T_k)) & t \in [0, T_k - T_{k-1}] \\ x_{k-1}(t - (T_k - T_{k-1}), x_k(T_k - T_{k-1}, z(T_k))) & t \in [T_k - T_{k-1}, T_k - T_{k-2}] \\ \vdots & \vdots \\ x_1(t - (T_k - T_1), x_2(T_2 - T_1, \ldots x_k(T_k - T_{k-1}, z(T_k)))) & t \in [T_k - T_1, T_k] \end{cases}$$

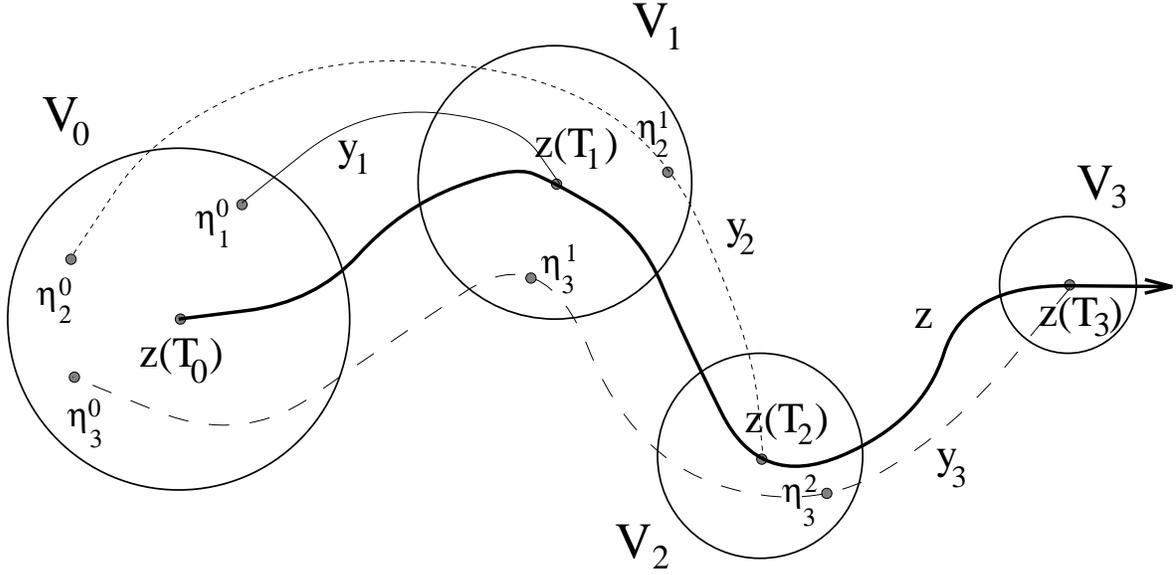

Figure 1: Construction of $\eta_j^k$

We set

$$\eta_0^0 = z(T_0)$$
$$\eta_1^0 = y_1(T_1),\ \eta_1^1 = z(T_1)$$
$$\eta_2^0 = y_2(T_2),\ \eta_2^1 = y_2(T_2 - T_1),\ \eta_2^2 = z(T_2)$$
$$\eta_3^0 = y_3(T_3),\ \eta_3^1 = y_3(T_3 - T_1),\ \eta_3^2 = y_3(T_3 - T_2),\ \eta_3^3 = z(T_3)$$
$$\vdots$$
$$\eta_j^0 = y_j(T_j),\ \ldots, \eta_j^k = y_j(T_j - T_k),\ \ldots,\ \eta_j^j = z(T_j)$$
$$\vdots$$

By construction, each $\eta_j^k \in V_k$ (see Figure 1).

It remains to verify that this construction satisfies the final condition. Suppose that for some $k \geq 1$, the subsequence $\{\eta_{j_l}^k\}_{l=1}^\infty$ converges to a limit $\eta^k$. Recall that by definition,

$$\eta_j^{k-1} = x_k(T_k - T_{k-1}, \eta_j^k)$$

for each $j \geq 0$. Then, by continuity of $x_k(T_k - T_{k-1}, \cdot)$, we find that $\{\eta_{j_l}^{k-1}\}_{l=1}^\infty$ is convergent, since

$$\lim_{l \to \infty} \eta_{j_l}^{k-1} = \lim_{l \to \infty} x_k(T_k - T_{k-1}, \eta_{j_l}^k) = x_k(T_k - T_{k-1}, \eta^k).$$

Denote $\eta^{k-1} = x_k(T_k - T_{k-1}, \eta^k)$. Finally, we note that $\overline{x}(t) : [0, T_k - T_{k-1}] \to X$ defined by $\overline{x}(t) := x_k(T_k - T_{k-1} - t, \eta^k)$ is a solution of (11), satisfies (12), and has $x(T_k - T_{k-1}) = \eta^k$. All that remains is to re-number the sequences $\eta_j^k$ so they each begin at $j = 1$. ∎

We now turn to our main result. To apply Lemma 3.2, we restrict to a setting in which the constructed sequences are guaranteed to have convergent subsequences; we suppose $X$ is finite dimensional. For a given $0 < T \leq \infty$ and a reference trajectory $z : [0, T) \to X$ of the relaxed system (2), we will construct a trajectory of the original system which stays within a given tube (with possibly vanishing radius) around the trajectory $z$.

**Theorem 1** *Suppose the space $X$ is finite dimensional. Let $0 < T \leq \infty$. Suppose the set-valued map $F : [0, T) \times X \to \mathcal{P}(X)$ satisfies the hypotheses (H1'')-(H3'') of Lemma 3.2. Fix $\xi \in X$ and let $z : [0, T) \to X$ be a solution of*

$$\dot{x} \in \text{clco}\, F(t, x), \quad x(0) = \xi.$$

*Let $r : [0, T) \to \mathbb{R}$ be a continuous function satisfying $r(t) > 0$ for all $t \in [0, T)$. Then there exists an $\eta^0 \in B(\xi, r(0))$ and a solution $x : [0, T) \to X$ of*

$$\dot{x} \in F(t, x), \quad x(0) = \eta^0, \qquad (15)$$

*which satisfies*

$$|z(t) - x(t)| \leq r(t) \qquad \forall t \in [0, T).$$

*Proof.* Choose a strictly increasing sequence of times $T_k$ so that $T_0 = 0$ and $T_k \to T$. Let the sequence $\{\delta_k\}_{k=0}^\infty$ and, for each nonnegative integer $k$, the sequence $\{\eta_j^k\}_{j=0}^\infty$ be as in Lemma 3.2 for $F$, $z$, and $r(\cdot)$. For each $k \geq 0$, set $V_k = B(z(T_k), \delta_k)$ and $r_{k+1} = \min\{r(t) : t \in [T_k, T_{k+1}]\}$.

Since the sequence $\{\eta_j^0\}_{j=1}^\infty$ lies in $V_0$, which is compact, there is a convergent subsequence $\{\eta_{j_{l0}}^0\}_{l0=1}^\infty$ which converges to some $\eta^0 \in V_0$. Likewise, the sequence $\{\eta_{j_{l0}}^1\}_{l0=1}^\infty$ lies in the compact set $V_1$, so it has a subsequence $\{\eta_{j_{l1}}^1\}_{l1=1}^\infty$ converging to some $\eta^1 \in V_1$. Continuing with this diagonalization, we find, for each $k \geq 1$, a subsequence $\{\eta_{j_{lk}}^k\}_{lk=1}^\infty$ which converges to some $\eta^k \in V_k$ and which is a subsequence of $\{\eta_{j_{l(k-1)}}^{k-1}\}_{l(k-1)=1}^\infty$.

Now, for each $k \geq 1$, since $\{\eta_{j_{lk}}^k\}_{lk=1}^\infty$ converges to $\eta^k$, it follows from Lemma 3.2 that there is a trajectory $\overline{x}_k : [0, T_k - T_{k-1}] \to X$ which solves (11), satisfies (12), and has $\overline{x}_k(T_k - T_{k-1}) = \eta^k$.

We construct a trajectory $x : [0, T) \to X$ by

$$x(t) := \overline{x}_k(t) \quad \text{when} \quad t \in [T_{k-1}, T_k).$$

By construction, this trajectory is a solution of (15) on the interval $[0, T)$ and satisfies $x(0) = \eta^0$. It follows from property (12) that

$$|z(t) - x(t)| \leq r_k \qquad \forall t \in [T_{k-1}, T_k],$$

from which we conclude

$$|z(t) - x(t)| \leq r(t) \qquad \forall t \in [0, T).$$

∎

Before stating a corollary, we quote a standard existence result for compact valued differential inclusions which follows from, e.g. [2], Theorem 2.3.1.

**Lemma 3.3** *Let $\Omega \subset \mathbb{R} \times \mathbb{R}^n$ be an open set containing $(0, x_0)$, and let $F$ be a locally Lipschitz set-valued map from $\Omega$ to the nonempty compact subsets of $\mathbb{R}^n$. Then there exists $\tau > 0$ and a solution $x$ of*

$$\dot{x} \in F(t, x), \quad x(0) = x_0$$

*defined on the interval $[0, \tau)$.* □

**Corollary 3.4** *Suppose $X = \mathbb{R}^n$ and the set-valued map $F : [0, \infty) \times \mathbb{R}^n \to \mathcal{P}(\mathbb{R}^n)$ has compact values and satisfies the hypotheses (H1'')-(H3'') of Lemma 3.2. Then the inclusion (1) is forward complete if and only if its relaxation (2) is forward complete.*

*Proof.* One implication is immediate. Suppose now that the inclusion (1) is forward complete but its relaxation (2) is not. Choose a maximal solution $z$ of (2) which has a bounded interval of definition $[0, T)$. Applying Theorem 1 with $r(t) = T - t$, we choose a solution $y$ of (1) on $[0, T)$ which satisfies

$$|y(t) - z(t)| \leq T - t \qquad \forall t \in [0, T). \tag{16}$$

Now, since the inclusion (1) is forward complete, the solution $y$ has an extension to the interval $[0, \infty)$, which we also call $y$. Since $y$ is continuous at $T$, we have, from (16),

$$\lim_{t \to T^-} z(t) = y(T).$$

By Lemma 3.3, there exists $\tau > 0$ so that

$$\dot{x} \in \operatorname{clco} F(T + t, x), \qquad x(0) = y(T)$$

has a solution $\widehat{z}$ defined on $[0, \tau)$. The concatenation of $z$ with $\widehat{z}$ is an extension of $z$ to the interval $[0, T+\tau)$ which contradicts the definition of $T$. We conclude that each maximal solution of (2) is defined on $[0, \infty)$, that is, (2) is forward complete. ∎

## 4 Counter-example

Considering that the Theorem in this note provides a complementary result to the classical Filippov-Ważewski Theorem, it is natural to ask whether one can achieve the results of both theorems simultaneously. That is, whether there exists an infinite-time approximation which satisfies the same initial condition as a given reference trajectory. The following example shows that in general this is not possible.

Consider the following differential inclusion evolving on $\mathbb{R}^2$:

$$\begin{aligned} \dot{x}(t) &= y^2(t) \\ \dot{y}(t) &\in \{-1, 1\}, \end{aligned}$$

and the relaxation to convex values:

$$\begin{aligned} \dot{x}(t) &= y^2(t) \\ \dot{y}(t) &\in [-1, 1]. \end{aligned}$$

Note that $x(t) \equiv y(t) \equiv 0$ is a solution of the relaxed inclusion with $x(0) = y(0) = 0$.

Clearly, the set-valued function $F(x, y) = (\{y^2\}, \{-1, 1\})$ is measurable, locally bounded, locally Lipschitz, and has closed, nonempty values. Then, by Theorem 1, the original inclusion admits solutions which approximate the zero solution for $t \geq 0$. For example, there exists a solution $(x(t), y(t))$ which satisfies

$$x^2(t) + y^2(t) \leq e^{-t} \qquad \forall t \geq 0, \tag{17}$$

with $|x(0)| \leq 1$, $|y(0)| \leq 1$.

However, the inclusion *cannot* admit a solution satisfying (17) and also satisfying $x(0) = y(0) = 0$. We note that any solution with $x(0) = y(0) = 0$ satisfies

$$x(1) = \int_0^1 y^2(t) \, dt = \varepsilon > 0,$$

for some $\varepsilon > 0$, as $\int_0^1 y^2(t) \, dt = 0$ implies $y(t) = 0$ almost everywhere on $[0, 1]$, which is not allowed. Then as $\dot{x}(t) \geq 0$ for all $t \geq 0$, it follows that $x(t) \geq \varepsilon$ for all $t \geq 1$, so (17) cannot be achieved.

## A  Appendix

For completeness, we state the main result in [4], which is the primary tool used in the proof of Theorem 1.

**Definition A.1** A set-valued map $F$ from a metric space $Z$ to subsets of a metric space $Y$ is called *lower semicontinuous* at $z \in Z$ if $F(z) \neq \emptyset$ and for any $y \in F(z)$ and any neighbourhood $N(y)$ of $y$, there exists a neighbourhood $N(z)$ of $z$ so that

$$F(\zeta) \cap N(y) \neq \emptyset \qquad \forall \zeta \in N(z).$$

The map $F$ is called lower semicontinuous if it is lower semicontinuous at each $z \in Z$.

Let $S$ be a separable metric space. Let $F : [0, 1] \times X \times S \to \mathcal{P}(X)$ and consider the following initial value problems

$$\dot{x} \in F(t, x, s), \qquad x(0) = \xi(s), \qquad (18)$$

where $\xi : S \to X$ is a continuous function.

**Theorem 2** Suppose the set-valued map $F$ satisfies

(H1) $F$ is $\mathcal{L}[0,1] \otimes \mathcal{B}(X \times S)$ measurable;

(H2) for any $(t, x)$, the map $s \mapsto F(t, x, s)$ is lower semicontinuous;

(H3) there exists a map $s \mapsto k(\cdot, s)$ continuous from $S$ into $L^1([0, 1], \mathbb{R})$ such that for any $s \in S$ and $\xi, \eta \in X$,

$$d_H(F(t, \xi, s), F(t, \eta, s)) \leq k(t, s) |\xi - \eta| \qquad \text{a.e. } t \in [0, 1];$$

(H4) for any continuous map $s \mapsto y(\cdot, s)$ from $S$ into $AC([0, 1], X)$, there exists a continuous map $\beta_y : S \to L^1([0, 1], \mathbb{R})$ such that for any $s \in S$,

$$d(\dot{y}(t, s), F(t, y(t, s), s)) \leq \beta_y(s)(t) \qquad \text{a.e. } t \in [0, 1]. \qquad (19)$$

Then for any continuous map $s \mapsto y(\cdot, s)$ from $S$ into $AC([0, 1], X)$, any map $s \mapsto \beta(s) = \beta_y(s)$ from $S$ into $L^1([0, 1], \mathbb{R})$ satisfying (19), and any $\varepsilon > 0$, there exists a function $x : [0, 1] \times S \to X$ such that

(a) for every $s \in S$, the function $t \mapsto x(t, s)$ is a solution of (18);

(b) the map $s \mapsto x(\cdot, s)$ is continuous from $S$ into $AC([0, 1], X)$;

(c) for every $s \in S$, and almost every $t \in [0, 1]$,
$$|\dot{y}(t,s) - \dot{x}(t,s)| \leq \varepsilon + \varepsilon k(t,s)e^{m(t,s)} + k(t,s)\,|y(0,s) - \xi(s)|\, e^{m(t,s)}$$
$$+ k(t,s) \int_0^t \beta(s)(\tau) e^{m(t,s) - m(\tau,s)}\, d\tau + \beta(s)(t)$$

(d) for every $s \in S$ and every $t \in [0, 1]$,
$$|[y(t,s) - x(t,s)] - [y(0,s) - \xi(s)]| \leq \varepsilon e^{m(t,s)} + |y(0,s) - \xi(s)|\,(e^{m(t,s)} - 1)$$
$$+ \int_0^t \beta(s)(\tau) e^{m(t,s) - m(\tau,s)}\, d\tau,$$

where $m(t, s) := \int_0^t k(\tau, s)\, d\tau$. □

**Remark A.2** By assumption (H3), the assumption (H4) can be replaced by the equivalent condition:
(H4$_0$) there exists a continuous map $\beta_0 : S \to L^1([0,1], \mathbb{R})$ such that for any $s \in S$
$$d(0, F(t, 0, s)) \leq \beta_0(s)(t) \quad \text{a.e. } t \in [0, 1].$$

**Remark A.3** The dependence of $F$ on the parameter $s$ is dropped for the purposes of this note, since the proof of Theorem 1 will not hold for this more general case. In that proof, the approximating trajectory is constructed by concatenating solutions which correspond to different values of $s$.

The authors would like to thank Héctor Sussmann for bringing the crucial references to our attention, and also David Angeli for helpful discussions.

# References

[1] D. Angeli, B. Ingalls, E. D. Sontag, and Y. Wang, *A Relaxation Theorem for Asymptotically Stable Differential Inclusions*, in preparation.

[2] J.-P. Aubin and A. Cellina, *Differential Inclusions*, Spring-Verlag, Berlin, 1984.

[3] J.-P. Aubin and H. Frankowska, *Set-Valued Analysis*, Birkhäuser, Boston, 1990.

[4] R. M. Colombo, A. Fryszkowski, T. Rzeżuchowski, and V. Staicu, *Continuous Selections of Solution Sets of Lipschitzean Differential Inclusions*, Funkcialaj Ekvacioj, **34** (1991), pp. 321–330.

[5] K. Deimling, *Multivalued Differential Equations*, Walter De Gruyter & Co., Berlin, 1992.

[6] A. F. Filippov, *Differential Equations with Discontinuous Righthand Sides*, Kluwer Academic, Dordrecht, The Netherlands, 1988.

[7] A. Fryszkowski and T. Rzeżuchowski, *Continuous Version of Filippov-Ważewski Relaxation Theorem*, Journal of Differential Equations, **94** (1991), pp. 254–265.

[8] E. D. Sontag and Y. Wang, *New characterizations of the input to state stability property*, IEEE Transactions on Automatic Control, **41** (1996), pp. 1283–1294.